\documentclass[11pt]{amsart}
\usepackage{times,epsf,amssymb,amsmath,hyperref, tikz, rlepsf, color, pstricks}

\begin{document}

\newtheorem{thm}{Theorem}[section]
\newtheorem{lem}[thm]{Lemma}
\newtheorem{cor}[thm]{Corollary}
\newtheorem{conj}[thm]{Conjecture}
\newtheorem{question}[thm]{Question}

\theoremstyle{definition}
\newtheorem{defn}[thm]{\bf{Definition}}

\theoremstyle{remark}
\newtheorem{rmk}[thm]{Remark}

\def\square{\hfill${\vcenter{\vbox{\hrule height.4pt \hbox{\vrule width.4pt height7pt \kern7pt \vrule width.4pt} \hrule height.4pt}}}$}

\newenvironment{pf}{{\it Proof:}\quad}{\square \vskip 12pt}

\title{Existence of Embedded Minimal Disks}
\author{Baris Coskunuzer}
\address{UT Dallas, Dept. Math. Sciences, Richardson, TX 75080}
\email{baris.coskunuzer@utdallas.edu}
\thanks{The author is partially supported by Simons Collaboration Grant, and Royal Society Newton Mobility Grant.}

\maketitle


\newcommand{\BHH}{\mathbf{H}^3}
\newcommand{\CH}{\mathcal{C}(\Gamma)}
\newcommand{\BH}{\mathbf{H}}
\newcommand{\BR}{\mathbf{R}}
\newcommand{\BC}{\mathbf{C}}
\newcommand{\BZ}{\mathbf{Z}}

\newcommand{\e}{\epsilon}

\newcommand{\wh}{\widehat}
\newcommand{\wt}{\widetilde}

\newcommand{\A}{\mathcal{A}}
\newcommand{\T}{\mathcal{T}}
\newcommand{\C}{\mathcal{C}}
\newcommand{\D}{\mathcal{D}}
\newcommand{\X}{\mathcal{X}}
\newcommand{\Y}{\mathcal{Y}}
\newcommand{\Z}{\mathcal{Z}}
\newcommand{\I}{\mathcal{I}}
\newcommand{\B}{\mathbf{B}}
\newcommand{\h}{\mathbf{h}}
\newcommand{\kk}{\mathbf{k}}

\begin{abstract}

We give a generalization of Meeks-Yau's celebrated embeddedness result for the solutions of the Plateau problem for extreme curves.

\end{abstract}

\section{Introduction}

The Plateau problem asks the existence of a least area disk for a given curve in the ambient manifold $M$. This problem was solved for $\BR^3$ by Douglas \cite{Do}, and Rado \cite{Ra} in early 1930s. Later, it was generalized by Morrey \cite{Mo} for Riemannian manifolds. After several results on the regularity of the solutions, embeddedness question was studied by many experts: {\em For which curves, the solution to the Plateau problem is embedded?} (See Figure \ref{nonembeddedfig}). 

After several partial results in 1970s \cite{GS,TT,AS}, Meeks-Yau \cite{MY1} finally showed that for extreme curves, the solution to the Plateau problem must be embedded. In 2002, Ekholm, White, and Wienholtz prove that for the curves of total curvature less than $4\pi$, the least area disks are embedded, too \cite{EWW}.

In this paper, we give a generalization of Meeks-Yau's celebrated embeddedness result for extreme curves as follows:

\begin{thm} Let $\Gamma$ be a weak-extreme curve in $\BR^3$. Then, $\Gamma$ bounds an embedded stable minimal disk $\Sigma$ in $\BR^3$.	
\end{thm}

A weak-extreme curve is basically a curve which naturally decomposes into extreme pieces. For formal definition, see Section \ref{weakextremesec}.

Furthermore, we discuss other potential directions, and generalizations of the embeddedness question for the solutions of the Plateau problem. In particular, we show that for the curves with convex hull genus $0$, the solution to the Plateau problem may not be embedded. On the other hand, while our main result shows that weak-extreme curves bound stable embedded minimal disks, the solution to the Plateau problem (least area disk) may not be embedded. We construct explicit counterexamples for both cases. By giving this generalization, and discussing these counterexamples, we aim to bring attention to this long-standing problem again.

The organization of the paper is as follows. In the next section, we give basic definitions, and the background. In Section 3, we prove our main result. In Section 4, we  construct our counterexamples mentioned above. In Section 5, we  give concluding remarks about further directions, and questions.

\subsection{Acknowledgements:}
We are very grateful to the referee for very valuable comments and suggestions. 

\section{Preliminaries}

In this section, we give the basic definitions, and overview the necessary background on the problem. For further details, see \cite{HS}.

\begin{defn}  A {\em least area disk} is a disk which has the smallest area among the disks with the same boundary.
An {\em area minimizing surface} is a surface which has the smallest area among all orientable surfaces
(with no topological restriction) with the same boundary. We call a surface {\em minimal} if the mean curvature vanishes everywhere.
\end{defn}

\begin{defn} [Convex Hull] Let $\Gamma$ be a curve (or arc) in $\BR^3$. Let $CH(\Gamma)$ be the smallest convex set in $\BR^3$ containing $\Gamma$. We call $CH(\Gamma)$ {\em convex hull of $\Gamma$}. 
\end{defn}

\begin{defn} [Mean Convexity] Let $\Omega$ be a compact region in $\BR^3$ with boundary. We call $\Omega$ a {\em mean convex domain} if the following conditions hold.

\begin{itemize}

\item $\partial \Omega$ is piecewise smooth.

\item Each smooth subsurface of $\partial \Omega$ has nonnegative mean curvature with respect to inward normal.

\item Each smooth subsurface $S$ of $\partial \Omega$  extends to a smooth embedded surface $\wh{S}$ in $\BR^3$ such
that $\wh{S} \cap \Omega = S$.

\end{itemize}

\end{defn}

Note that mean convexity definition is indeed more general, and applies to general Riemannian $3$-manifolds \cite{MY2,HS}. We just adapted this simpler definition as we mainly work in $\BR^3$.

\begin{figure}[t]
	
	\relabelbox  {\epsfxsize=3in
		
		\centerline{\epsfbox{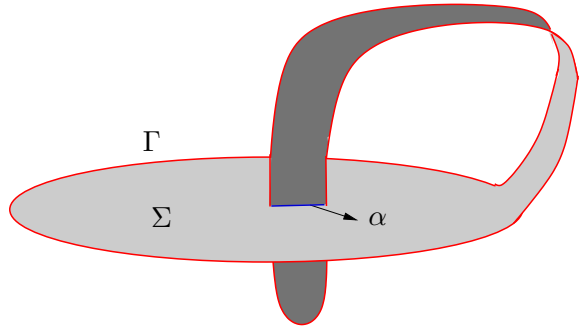}}}

	\relabel{1}{ $\Gamma$} 
	\relabel{2}{ $\Sigma$} 
	\relabel{3}{$\alpha$}

	\endrelabelbox
	
	\caption{\label{nonembeddedfig} \small The red curve $\Gamma$ is an example of a curve, which is not extreme. Furthermore, the least area disk $\Sigma$ with $\partial \Sigma=\Gamma$ has self-intersection along $\alpha$.}
	
\end{figure}

\begin{defn} [Extreme Curve] A simple closed curve $\Gamma$ is an {\em extreme curve} if it is on the boundary of its convex hull $CH(\Gamma)$, i.e. $\Gamma\subset \partial CH(\Gamma)$. We  call an arc $\gamma$ \textit{extreme} if $\gamma\subset \partial CH(\gamma)$.
\end{defn}

Now, we state the former results which we use in the following sections. While their result applies to general Riemannian $3$-manifolds, we give their simpler versions adapted to $\BR^3$.

\begin{lem}\cite{MY1}, \cite{MY2} \label{MYlem}
Let $\Omega$ be a compact, mean convex domain in $\BR^3$. Let $\Gamma\subset\partial \Omega$ be a  simple closed curve which is nullhomotopic in $\Omega$. Then, there exists a least area disk $D\subset \Omega$ with $\partial  D = \Gamma$. Moreover, all such disks are
properly embedded in $\Omega$ and they are pairwise disjoint. 
\end{lem}

\begin{rmk} The above lemma implies the solutions to the Plateau problem for extreme curves are embedded. In particular, any minimal surface bounding a curve $\Gamma$ must belong to the convex hull $CH(\Gamma)$ by maximum principle. Hence, the least area disk in the convex hull $CH(\Gamma)$ will be the least area disk in $\BR^3$, too. By the lemma above, the least area disk in the convex hull must be embedded. See Figure \ref{nonembeddedfig} to see why the extremeness condition is crucial for embeddedness. 
	
Note that only mean convexity condition is not enough to secure the embeddedness of the least area disk in the ambient manifold.  In \cite{Co1}, we constructed examples of curves in the boundary of a mean convex domain in $\BR^3$ where the solution to the Plateau problem in $\BR^3$ is not embedded.	See also Remark \ref{rmkLAdisk}.
\end{rmk} 

Furthermore, Meeks-Yau proved a partial converse of this statement as follows:

\begin{lem}\cite{MY3}
If a simple closed curve $\Gamma$ bounds a strictly stable embedded minimal disk, then there exists a mean convex $3$-manifold $M$ such that $\Gamma\subset \partial M$.
\end{lem}

Moreover, there is an analogous existence result for area minimizing surfaces.

\begin{lem} \label{AMSlem} \cite{Fe}, \cite{HSi}, \cite{Wh} 	
Let $\Omega$ be a compact, mean convex domain in $\BR^3$. Let $\Gamma$ be a  nullhomologous, simple closed curve  in $\Omega$. Then, there exists a smoothly embedded area minimizing surface $\Sigma\subset \Omega$ with $\partial \Sigma =\Gamma$.
\end{lem}

\section{Existence of Embedded Stable Minimal Disks for Weak-Extreme Curves}

In this section, we  first define the weak-extreme curves which is a generalization of extreme curves. Then, we show that they bound stable embedded minimal disks.

\subsection{Weak-Extreme Curves} \label{weakextremesec} \

\vspace{.2cm}

Let $\Gamma$ be a simple closed curve in $\BR^3$. Let $\X=CH(\Gamma)$ be the  convex hull of $\Gamma$. Let $\Gamma^*=\Gamma\cap \partial \X$. If $\Gamma^*$ have finitely many components, then we call $\Gamma$ {\em tame}, otherwise {\em wild}. Throughout the paper, we  only deal with {\em tame curves}. So, assume $\Gamma$ is tame. 

Now, let $\Gamma^*=\alpha_1\cup\alpha_2 .. \cup\alpha_n$ where  $\{\alpha_i\}$ are the connected components of $\Gamma^*$. Here, $\alpha_i$ is either an arc or an isolated point. If $\alpha_i$ is an arc, let $\partial \alpha_i=\{p_{2i-1},p_{2i}\}$. In the case $\alpha_i$ is just a point $q$, let $p_{2i-1}=p_{2i}=q$.

Let $\Gamma-\Gamma^*=\beta_1'\cup\beta_2'\cup ... \cup\beta_n'$ where $\beta_i$ is the connected open arc in $\Gamma$ such that  $\partial \overline{\beta_i'}=\{p_{2i},p_{2i+1}\}$ (use mod $2n$ for $i=n$). Let $\beta_i=\overline{\beta_i'}$. 



\begin{figure}[b]
	\begin{center}
		$\begin{array}{c@{\hspace{.2in}}c}

		\relabelbox  {\epsfxsize=2.1in \epsfbox{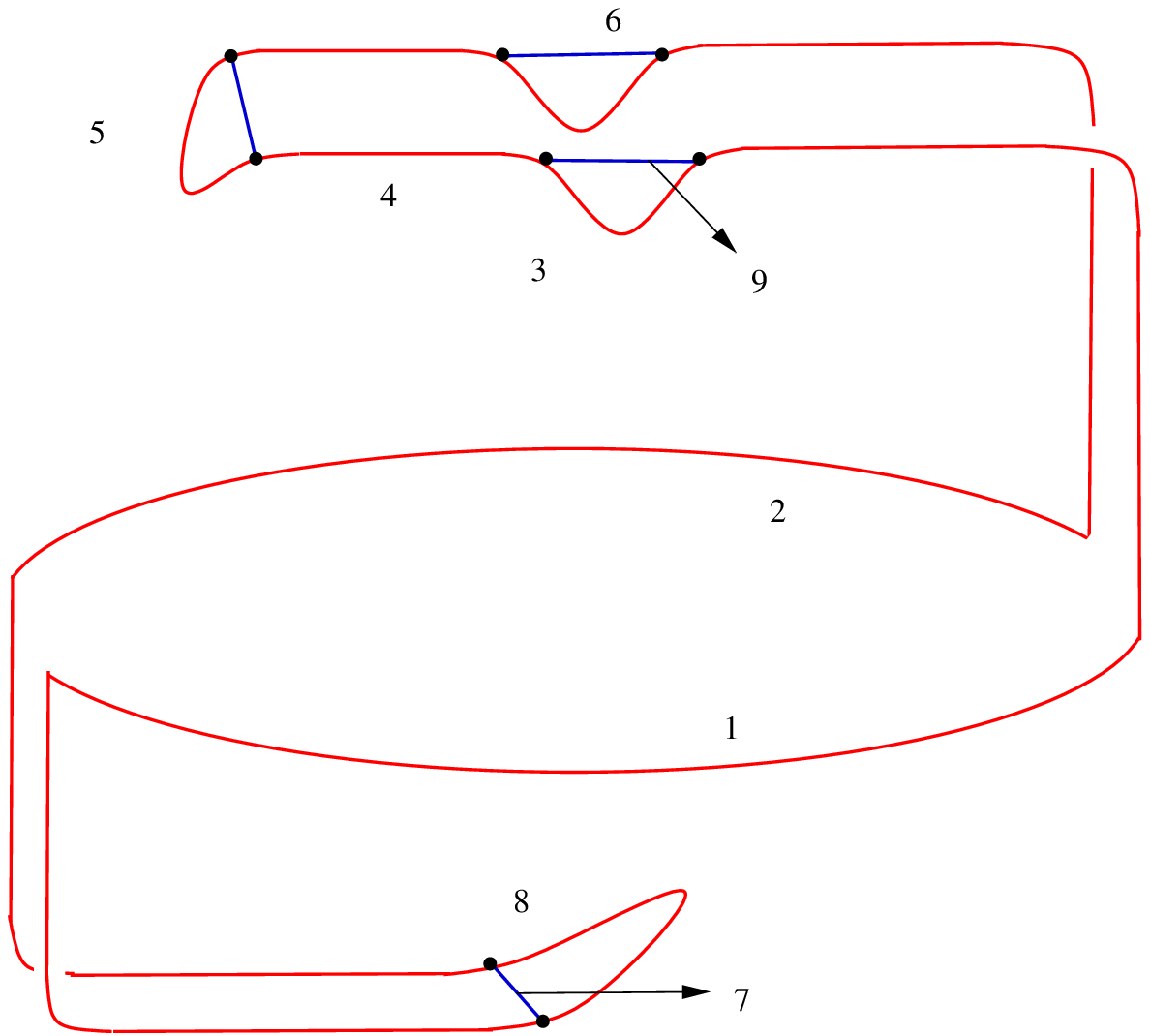}}
		\relabel{1}{\footnotesize $\alpha_1$} 
		\relabel{2}{\footnotesize $\alpha_4$}  
		\relabel{3}{\footnotesize $\beta_1$} 
		\relabel{4}{\footnotesize $\alpha_2$} 
		\relabel{5}{\footnotesize $\beta_2$} 
		\relabel{6}{\footnotesize $l_3$} 	
		\relabel{7}{\footnotesize $l_4$} 	   
		\relabel{8}{\footnotesize $\beta_4$} 	
		\relabel{9}{\footnotesize $l_1$} 	
		\endrelabelbox &
		
		\relabelbox  {\epsfxsize=1.5in \epsfbox{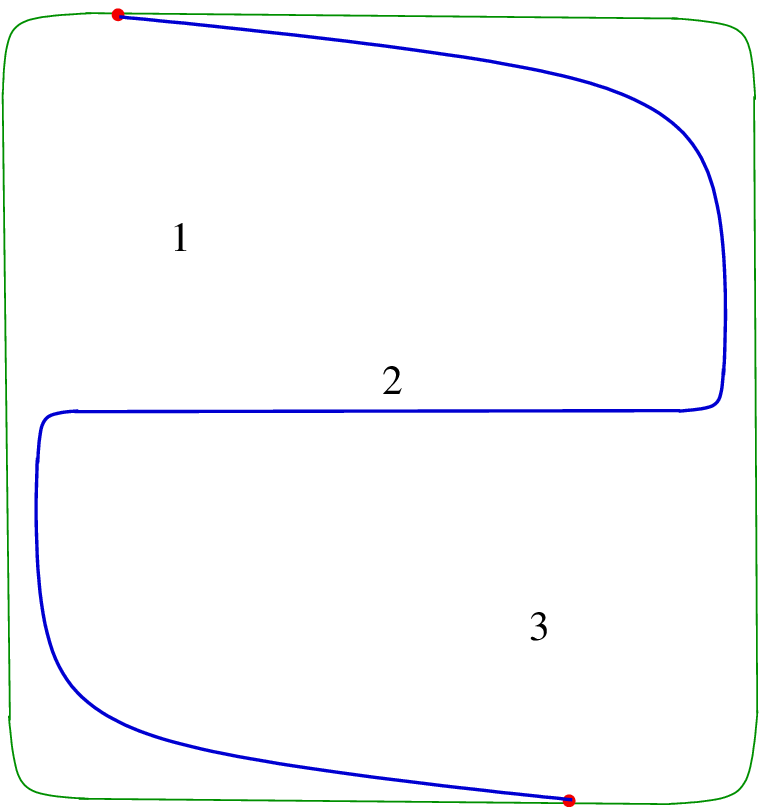}}
		\relabel{1}{$\Omega^+$} 
		\relabel{2}{$\wh{\Sigma}$}  
		\relabel{3}{$\Omega^-$} 
		\endrelabelbox \\
		
	\end{array}$
	
\end{center}

\caption{ \label{fig_weakextreme} \footnotesize In the figure left, we decompose red curve $\Gamma$ into extreme pieces $\alpha_1\cup\beta_1\cup...\cup\alpha_4\cup\beta_4$. Here, the hooks $\beta_1, \beta_2, \beta_3$ belongs to $\Omega^+$, while $\beta_4$ belongs to $\Omega^-$. In the figure right, we depict convex hull $CH(\Gamma)$ in 2D way from the side. Green box represents $CH(\Gamma)$, while blue arc represents $\wh{\Sigma}$. }
\end{figure}

Now, we define an extreme Jordan curve $\wh{\Gamma}$ in $\partial \X$ by modifying $\Gamma$. Let $l_i$ be a path in $\partial \X$ with the same endpoints with $\beta_i$, i.e. $\partial \beta_i=\partial l_i=\{p_{2i},p_{2i+1}\}$, and $l_i\subset \partial \X$.  As $\Gamma^*$ is not a closed loop in $\partial \X$, we can always choose $\{l_i\}$ inductively so that $\wh{\Gamma}=\Gamma^*\bigcup_{i=1}^n l_i$ is embedded. Hence, this gives us an extreme simple closed curve $\wh{\Gamma}$ in $\partial \X$. See Figure \ref{fig_weakextreme}-left.  Note that whenever possible, we will  choose $l_i$ as the line segment connecting $p_{2i}$ and $p_{2i+1}$, which simplifies the construction. See Remarks \ref{rmk-gammahead}, \ref{rmk-l_i-notstraight}, \ref{rmk-alternative} for further discussion.

Now, for any $i$, define $\wh{\beta}_i=\beta_i\cup l_i$. As $l_i \cap \beta_i= \{p_{2i},p_{2i+1}\}$, $\wh{\beta}_i$ is a simple closed curve in $\X$. In our construction, we will assume $\wh{\beta}_i$ is an extreme curve, and call $\wh{\beta}_i$ {\em hooks} of $\Gamma$.

In particular, we obtain an {\em extremal decomposition} of $\Gamma$ as follows: We have $n+1$ extreme curves, i.e. $\wh{\Gamma}$, and $\wh{\beta}_i$ for $1\leq i \leq n$. We have $\wh{\Gamma}\cap \bigcup_{i=1}^n \wh{\beta}_i = \bigcup_{i=1}^n l_i$. Furthermore, $\Gamma$ is the closure of $\wh{\Gamma}\Delta \bigcup_{i=1}^n \wh{\beta}_i$ where $\Delta$ is the symmetric difference.

Now, by construction, $\wh{\Gamma}$ is an extreme simple closed curve in $\partial \X$ such that $CH(\wh{\Gamma})=CH(\Gamma)=\X$. Hence, by Lemma \ref{MYlem}, any least area disk $\wh{\Sigma}$ with $\partial \wh{\Sigma}=\wh{\Gamma}$ is embedded. We  call the least area disk $\wh{\Sigma}$ {\em the core of the $CH(\Gamma)$}. See Figure \ref{fig_weakextreme}-right.

Now, we are ready to give the definition of weak-extreme curves.

\begin{defn} [Weak Extreme Curves] \label{wextdefn}
Let $\Gamma,\wh{\Gamma},\wh{\Sigma}, \{\beta_i\}, \{\wh{\beta}_i\}$ be as described above. We will call $\Gamma$ {\em weak-extreme}, if the following conditions hold.

\begin{enumerate}
	\item For any $i$, $\wh{\beta}_i$ is an extreme curve.
	\item $CH(\beta_i) \cap \beta_j=\emptyset$ for any $i\neq j$.
	\item $\beta_i\cap\wh{\Sigma}=\emptyset$ for any $i$.	
\end{enumerate}
\end{defn}

One can consider the above conditions as follows: 

{\em Condition 1:} Hooks are tight (extreme).

{\em Condition 2:} Hooks do not link each other.

{\em Condition 3:} Hooks do not stab the core of the convex hull (See Figure \ref{nonembeddedfig}).

After the proof of the main result, we explain the necessity and importance of these conditions for our construction in Remark \ref{rmk-Wecond}. See also an alternative definition for weak-extreme curves in Remark \ref{rmk-alternative}.

\subsection{Existence of Embedded Minimal Disks} \

\vspace{.2cm}

Our main result is as follows:

\begin{thm} \label{mainthm} Let $\Gamma$ be a weak-extreme curve in $\BR^3$. Then, $\Gamma$ bounds an embedded stable minimal disk $\Sigma$ in $\BR^3$.
\end{thm}

\noindent {\em Outline of the Proof:} Roughly, we first decompose our curve $\Gamma$ into extreme pieces, i.e. $\wh{\Gamma}$ and $\{\wh{\beta}_i\}$. Then, we get an embedded least area disk for each piece such that $\wh{\Sigma}$ with $\partial \wh{\Sigma}=\wh{\Gamma}$, and $\D_i$ with $\partial \D_i=\wh{\beta}_i$. By combining all these disks, we get a piecewise smooth embedded minimal disk $\Sigma'=\wh{\Sigma}\bigcup_{i=1}^n\D_i$. Then, we construct a special mean convex neighborhood $\Z$ around this piecewise smooth disk $\Sigma'$ by using area minimizing surfaces in both sides. Then, we get our stable embedded minimal disk $\Sigma$ in $\Z$. In particular, our smooth embedded minimal disk $\Sigma$ will basically be a "smoothening" of piecewise smooth disk $\Sigma'$. 	

\vspace{.2cm}

\begin{pf} Let $\Gamma, \wh{\Gamma},\Gamma^*, \wh{\Sigma}, \alpha_i,\beta_i,\wh{\beta}_i, l_i,\X$ be as defined above. As the convex hull $\X=CH(\Gamma)$ is a convex set in $\BR^3$, $\X$ is topologically a closed $3$-ball. Let $\Omega^\pm$ be the closure of the components of $\X-\wh{\Sigma}$, i.e. $\Omega^+\cap\Omega^-=\wh{\Sigma}$, and $\X=\Omega^+\cup\Omega^-$. Notice that both $\Omega^+$ and $\Omega^-$ are mean convex domains. See Figure \ref{fig_weakextreme}-right.
	
\vspace{.2cm}
	
\noindent {\bf Step 1:} {\em Construction of area minimizing surfaces $\T^\pm$ in $\Omega^\pm$:}

\vspace{.2cm} 

Recall that, by assumption, for any $i$, $\beta_i\cap\wh{\Sigma}=\emptyset$. Let $\I^+=\{i\mid \beta_i\subset\Omega^+\}$  and $\I^-=\{i\mid \beta_i\subset\Omega^-\}$, i.e. $\I^+\cup\I^-=\{1,2,...,n\}$.

Define $\wh{\Gamma}^+=\Gamma^*\bigcup_{i\in \I^+}\beta_i\bigcup_{j\in\I^-}\l_j$. Similarly, define $\wh{\Gamma}^-=\Gamma^*\bigcup_{i\in \I^-}\beta_i\bigcup_{j\in\I^+}\l_j$. In particular, we obtain $\wh{\Gamma}^+$ from $\wh{\Gamma}$ by replacing $l_i$ with $\beta_i$ for $i\in \I^+$. Similarly, we obtain $\wh{\Gamma}^-$ from $\wh{\Gamma}$ by replacing $l_j$ with $\beta_j$ for $j\in \I^-$. By construction, we have $\wh{\Gamma}^+\subset \Omega^+$ and $\wh{\Gamma}^-\subset \Omega^-$. 

Now, let $\T^+$ be the area minimizing surface in $\Omega^+$ with $\partial \T^+=\wh{\Gamma}^+$. Such an area minimizing surface exists by Lemma \ref{AMSlem}. Furthermore, $\T^+$ is smoothly embedded since $\Omega^+$ is mean convex. Similarly define the area minimizing surface $\T^-$ in $\Omega^-$ with $\partial \T^-=\wh{\Gamma}^-$.

\vspace{.2cm}

\noindent {\bf Step 2:} {\em Construction of least area disks $\D_i$ bounding hooks $\wh{\beta}_i$:} 

\vspace{.2cm}

Recall that $\wh{\beta}_i=\beta_i\cup l_i$. Let $\Y_i=CH(\wh{\beta_i})$. By assumption, $\wh{\beta}_i$ is an extreme curve, i.e. $\wh{\beta}_i\subset \partial \Y_i$. 

Let $i\in\I^+$. Consider $\Y_i-\T^+$. Let $\wh{\Y}_i$ be the component of $\Y_i-\T^+$ containing $l_i\subset \partial \X$. Notice that $\wh{\Y}_i$ is mean convex by construction. Let $\partial^+\Y_i$ and $\partial^-\Y_i$ be the components of $\partial \Y_i-\wh{\beta}_i$ where $\partial^-\Y_i$ is the component with $\partial^-\Y_i\cap \T^+\neq \emptyset$. We claim that $\partial^+\Y_i\cap \T^+=\emptyset$. This is because $\T^+$ is area minimizing surface in $\Omega^+$ and $l_i\subset \partial \Omega^+$. Hence, $l_i\cap\T^+=\emptyset$. On the other hand, as $\Y_i$  is convex, and $\Y_i\cap\beta_j=\emptyset$ for any $i\neq j$, the surface $\T^+$ cannot link $\beta_i$ from outside. In other words, $\partial \T^+\cap \Y_i =\beta_i$. Hence, 
$\T^+$ cannot have handles linking $\wh{\beta}_i$ as its projection to $\Y_i$ will have less area. This shows  $\partial^+\Y_i\cap \T^+=\emptyset$. 

Now, $\wh{\beta}_i$ is in $\partial \wh{\Y}_i$ and  nullhomotopic in $\wh{\Y}_i$ since $\partial^+ \Y_i$ is a disk with boundary $\wh{\beta}_i$. As $\wh{\Y}_i$ is mean convex, by Lemma \ref{MYlem}, $\wh{\beta}_i$ bounds a least area disk $\D_i$ in $\wh{\Y}_i$. By construction, $\D_i\cap \T^+=\emptyset$. Similarly, for any $j\in \I^-$, we get a least area disk in $\D_j$ in $\wh{\Y}_j$ with $\partial \D_j=\wh{\beta}_j$ such that $\D_j\cap \T^-=\emptyset$.

\vspace{.2cm}

\noindent {\bf Step 3:} {\em $\D_i\cap \D_j=\emptyset$ \ for any $i\neq j$.} 

\vspace{.2cm}

We prove this by using Meeks-Yau exchange roundoff trick. If $i\in \I^+$, and $j\in \I^-$, this is automatic as $\D_i\subset \Omega^+$ and $\D_j\subset \Omega^-$. Let $i,j\in \I^+$. Then, by assumption, $\wh{\beta}_i$ and $\wh{\beta}_j$ does not link each other. Furthermore, both of them are disjoint extreme curves. Hence, if $\D_i$ and $\D_j$ are the least area disks above, then $\D_i\cap \D_j$ is disjoint from their boundaries $\wh{\beta}_i$ and $\wh{\beta}_j$. Hence, if nonempty, $\D_i\cap \D_j$ contains a simple closed curve $\tau$. Let $E_i\subset \D_i$ and $E_j\subset \D_j$ with $\partial E_i=\partial E_j=\tau$. Consider $\D_i'=(\D_i-E_i)\cup E_j$. $\D_i$ and $\D_j$ being least area disks, we have $|E_i|=|E_j|$, and hence $|\D_i'|=|\D_i|$. This implies that $\D_i'$ is also a least area disk. However, $\D_i'$ has a folding curve along $\tau$. One can push $\D_i'$ toward convex side along the folding curve $\tau$, and reduce the area \cite{MY3}. This contradicts to $\D_i'$ being least area.  The claim follows.

\vspace{.2cm}

\noindent {\bf Step 4:} {\em Construction of mean convex domain $\Z$:} 

\vspace{.2cm}

Now, we are ready construct our mean convex domain $\Z$ in $\X$.  Recall that $\partial \T^\pm=\wh{\Gamma}^\pm$. Define a piecewise smooth surface $\wt{\T}^+=\T^+\bigcup_{i\in\I^+}\D_i$. As $\D_i\cap\D_j=\emptyset$ for any $i\neq j$, and $\T^+\cap\D_i=\beta_i$ for any $i\in\I^+$, $\wt{\T}^+$ is embedded. By construction, $\partial \wt{\T}^+=\wh{\Gamma}$. In particular, for each $i\in\I^+$, by adding $\D_i$ to $\T^+$, we "push" (replace) the boundary piece $\beta_i$ to $l_i$. Hence, after adding each $\D_i$ for $i\in \I^+$, we modify $\wh{\Gamma}^+=\Gamma^*\bigcup_{i\in \I^+}\beta_i\bigcup_{j\in\I^-}\l_j$ into $\wh{\Gamma}=\Gamma^*\bigcup_{i=1}^n l_i$. Similarly, define $\wt{\T}^-=\T^-\bigcup_{j\in\I^-}\D_j$. Again, $\partial \wt{\T}^-=\wh{\Gamma}$. 

Now, let $\wt{\T}=\wt{\T}^+\cup\wt{\T}^-$. As $\partial \wt{\T}^+=\partial \wt{\T}^-=\wh{\Gamma}$, $\wt{\T}$ is a closed embedded surface. As $\X$ is topologically a closed ball, any closed surface is separating in $\X$. Hence, define the region $\Z$ in $\X$ separated by $\wt{T}$, i.e. $\partial \Z=\wt{\T}=\T^+\cup\T^-\bigcup_{i=1}^n\D_i$. $\Z$ is mean convex as any smooth subsurface in $\partial \Z$ smoothly continues outside of $\Z$ by construction.

\vspace{.2cm}

\noindent {\bf Step 5:} {\em Existence of embedded, stable minimal disk $\Sigma$:} 

\vspace{.2cm}

Now, $\Gamma \subset \partial \Z$ by construction. Furthermore, $\Gamma$ is nullhomotopic in $\Z$ as $\Gamma$ is the boundary of the embedded (piecewise smooth) disk $\Sigma'=\wh{\Sigma}\bigcup_{i=1}^n \D_i$ in $\Z$. Therefore, by Lemma \ref{MYlem}, there exists a smoothly embedded least area disk $\Sigma$ in $\Z$. $\Sigma$ is a stable embedded minimal disk in $\BR^3$. The proof follows.	
\end{pf}

\begin{rmk} Note that it is possible to generalize this result to general Riemannian $3$-manifolds, by adapting convex hull, and definition of extremity  to this setting. We shall discussed this generalization in Section \ref{3mansec}.	
\end{rmk}

\subsection{Remarks on the Construction, and Definition of Weak-Extreme Curves} \label{sec-rmks}

\begin{rmk} \label{rmkLAdisk} {\em Being Least Area vs. Stable Minimal} 
	
\vspace{.2cm}
	
Note that the embedded minimal disk $\Sigma$ is a least area disk in $\Z$ by construction. However, it may not be a least area disk in $\BR^3$. We constructed examples of this type of curves in \cite{Co1}. In particular, while $\Gamma$ is in the boundary of a mean convex domain in $\BR^3$, the least area disk in $\BR^3$ bounding $\Gamma$ may be completely different, and it may not be embedded. In other words, being least area in a mean convex domain does not imply being least area in the larger space. In Section \ref{Weakcounterexamplesec}, we  discuss whether $\Gamma$ being weak-extreme implies the embeddedness of the solution to the Plateau problem in $\BR^3$.	
\end{rmk}

\begin{rmk} \label{rmk-Wecond} {\em The Necessity of the Conditions for Weak-Extreme Curves} 
	
\vspace{.2cm}
	
For our construction above, we need all 3 conditions in the definition of weak-extreme curves for the following reasons.
	
	\begin{enumerate}
		
		\item \textbf{Tightness of Hooks:} Our need $\wh{\beta}_i=\beta_i\cup l_i$ to be extreme is two fold. First reason is to make sure that the existence of least area disks $\D_i$ with $\partial \D_i=\wh{\beta}_i$. For example, if a  $\beta_i$ is too wiggly, then it night be impossible to have an embedded minimal disk bounding $\Gamma$. The second reason is that we need $\T^+\cap\D_i=\emptyset$ for our construction mean convex domain $\Z$. In Step 2, we crucially use extremeness of $\wh{\beta}_i$ to make sure that $\T^+\cap\D_i=\emptyset$.
		
		\item {\bf Hooks does not tangle each other:} $CH(\beta_i) \cap \beta_j=\emptyset$ for $i\neq j$. We needed this condition in Step 3, when we show $\D_i\cap \D_j=\emptyset$. Notice that even though hooks are disjoint, and $l_i\cap l_j=\emptyset$, the extreme hooks $\{\beta_i\}$ might tangle each other inside the convex hull $CH(\Gamma)$. This would make impossible to have $\D_i\cap\D_j=\emptyset$.
		
		
		\item {\bf Hooks do not stab the core:} $\beta_i\cap \wh{\Sigma}=\emptyset$. Considering Figure \ref{nonembeddedfig}, this is the most crucial condition for our construction. In that figure, the hook is intersecting the core $\wh{\Sigma}$, and hence the convex hull genus is no longer $0$. So, such a curve cannot bound an embedded minimal disk (See Section \ref{CHgenussec}).
		
	\end{enumerate}

\end{rmk}

\begin{rmk} \label{rmk-gammahead} {\em Dependence of $\wh{\Gamma}$ on the choice $\{l_i\}$} 

\vspace{.2cm}

Notice that in our definition for weak extreme curves, we have the flexibility in the choice of the extremal arcs $\{l_i\}$ in $\partial\X$ to replace $\{\beta_i\}$ when constructing $\wh{\Gamma}$. In these choices, there are two important conditions $\{l_i\}$ needs to satisfy. The first one is that $\beta_i\cup l_i$ needs to be extreme curve. Depending on $\beta_i$, being in $\partial X$ does not guarantee $\beta_i\cup l_i$ is extreme even though the arc $\beta_i$ is extreme. 
	
The second condition $\{l_i\}$ needs to satisfy is about the least area disk $\wh{\Sigma}$ which $\wh{\Gamma}$ bounds. Naturally, the disk depends on $\{l_i\}$ as $\wh{\Gamma}=\Gamma^*\bigcup_{i=1}^n l_i$. The third condition in our weak extreme curve definition requires $\beta_i\cap \wh{\Sigma}=\emptyset$. Hence, choice of $\{l_i\}$ is critical to make sure that $\wh{\Sigma}$ is away from the hooks.	
\end{rmk}

\begin{figure}[b]
	\begin{center}
		$\begin{array}{c@{\hspace{.2in}}c}

		\relabelbox  {\epsfxsize=2in \epsfbox{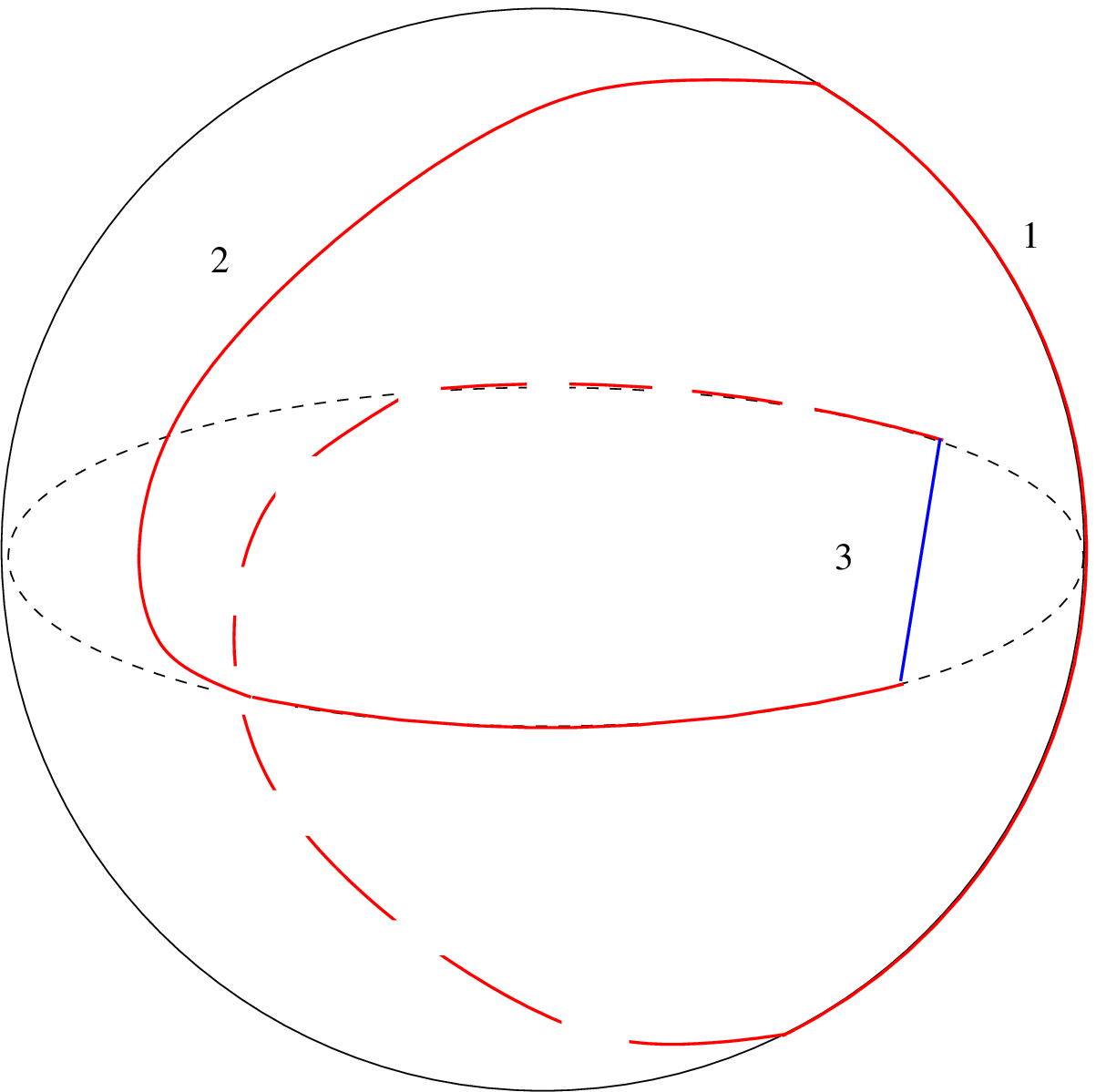}}
		\relabel{1}{\footnotesize $\Gamma$} 
		\relabel{2}{\footnotesize $\alpha_1$}  
		\relabel{3}{\footnotesize $\beta_1$} 
		
		\endrelabelbox &
		
		\relabelbox  {\epsfxsize=2in \epsfbox{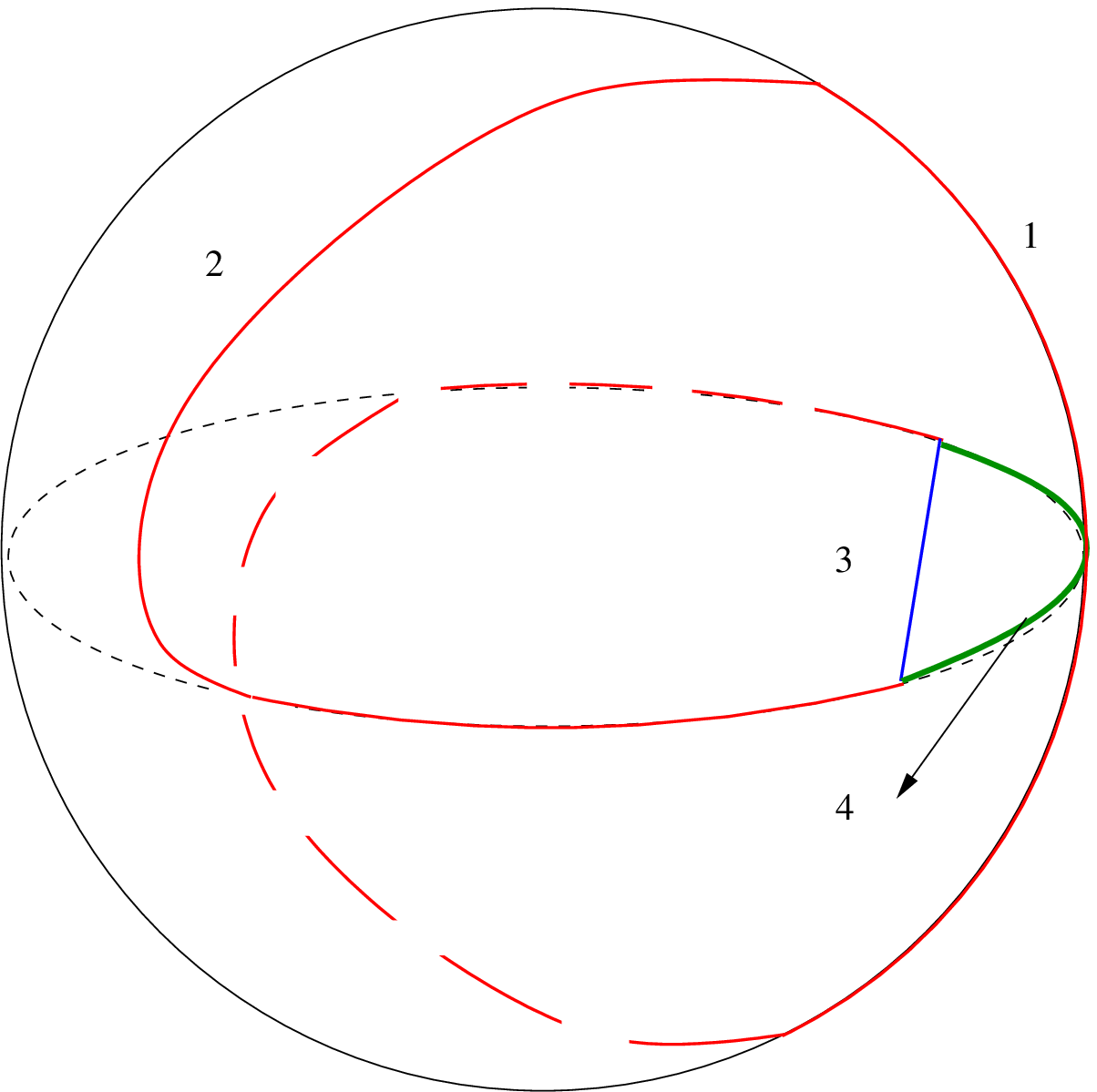}}
		\relabel{1}{\footnotesize $\Gamma$} 
		\relabel{2}{\footnotesize $\alpha_1$}  
		\relabel{3}{\footnotesize $\beta_1$} 
		\relabel{4}{\footnotesize $l_1$} 
		\endrelabelbox \\
		
	\end{array}$
	
\end{center}

\caption{ \label{fig_alternate} \footnotesize In the figure left, we decompose $\Gamma$ into two extreme pieces $\alpha_1\subset S^2$ (red arc) and $\beta_1$ (blue line segment). In the right, the green arc $l_1$ is the shortest path connecting $\partial\beta_1$ in $\partial \X$. However, the extreme curve $\wh{\Gamma}=\alpha_1\cup l_1$ is no longer embedded. }
\end{figure}

\begin{rmk} \label{rmk-l_i-notstraight} {\em $\{l_i\}$ not being straight line segments} 
	
\vspace{.2cm}
	
One good question at this point is that {\em Why don't we just choose $l_i$ as a straight line segment connecting the endpoints of corresponding $\beta_i$?} It is clear that such a choice would simplify the construction a lot, and it will make $\wh{\Gamma}$ and $\wh{\Sigma}$ canonical, only depending on $\Gamma$. In most cases, we can do that. However, when $\Gamma$ is a complicated curve (Figure \ref{fig_alternate}), the straight line segment $l_i$ may not be in the boundary of the convex hull $\partial \X$, but inside $\X$. This would imply that $\wh{\Gamma}$ is no longer an extreme curve, which is a highly crucial property for our construction. In Figure \ref{fig_alternate}, if we choose $l_1$ to be straight line segment connecting $\partial \beta_1$, then it would be the same line segment $l_1=\beta_1$, but $\wh{\Gamma}=\alpha_1\cup l_1$ would not be extreme anymore.
\end{rmk}

\begin{rmk} \label{rmk-alternative} {\em Alternative Construction: Canonical $\wh{\Sigma}$ while non-embedded $\wh{\Gamma}$}
	
\vspace{.2cm}
	
In our construction, $\wh{\Gamma}$, and $\wh{\Sigma}$ depend on our choices of extremal arcs $\{l_i\}$ in $\partial \X$ which are replacing the hooks $\{\beta_i\}$. There is an alternative construction which eliminates this freedom for choices of $\{l_i\}$. Hence, it gives us a canonical $\wh{\Gamma}$, and $\wh{\Sigma}$. For any $i$, we can require that the extreme arc $l_i$ to be the shortest path in $\partial \X$ connecting the endpoints of $\beta_i$. This requirement eliminates the randomness of $l_i$.

However, the problem with this requirement is that we might loose the embeddedness of $\wh{\Gamma}$ (See Figure \ref{fig_alternate}). Yet, this would not be a problem for our construction as embeddedness of $\wh{\Gamma}$ is not crucial for our construction. In this situation, one can modify our proof to get an embedded stable minimal disk as follows. 

If $\wh{\Gamma}$ is not embedded, then write $\wh{\Gamma}=\wh{\Gamma}_1\cup .. \cup \wh{\Gamma}_k$ where each $\wh{\Gamma}_j$ is a simple closed loop in $\partial \X$. Let $\wh{\Sigma}_i$ be the least area disk in $\X$ with $\partial \wh{\Sigma}_j=\wh{\Gamma}_j$. Then define $\wh{\Sigma}= \bigcup_{j=1}^k \wh{\Sigma}_j$. Here, $\wh{\Sigma}$ is a piecewise smooth disk with $\partial \wh{\Sigma}=\wh{\Gamma}$. 

Now, define $\X-\wh{\Sigma}=\Omega^+_1\cup .. \Omega^+_k\cup\Omega^-$, where we have $k+1$ components. Here, $\wh{\Sigma}_i\subset \partial \Omega^+_i$ and $\partial \Omega^-\supset \wh{\Sigma}$. In particular, instead of $\X-\wh{\Sigma}=\Omega^+\cup\Omega^-$ as in the original proof, we have a disconnected $\Omega^+$. By replacing the hooks $\beta_i$ in $\Omega^+_j$ with the corresponding $l_i$, we get simple closed curves $\wh{\Gamma}^+_j$ in $\Omega^+_j$. Then, we construct $k+1$ area minimizing surfaces $\T^+_j\subset \Omega_j^+$ with $\partial \T^+_j=\wh{\Gamma}^+_j$ for $1\leq j\leq k$ (Step 1 in Theorem \ref{mainthm}). Step 2 and Step 3 would be same.  For Step 4, we obtain mean convex domain $\Z$ in the same way with $\partial \Z =\T^-\cup \bigcup_{j=1}^k \T_j^+\bigcup_{i=1}^n\D_i$. Again, $\Gamma$ will be a nullhomotopic curve in $\partial \Z$. Then, we obtain a least area disk $\Sigma$ in $\Z$ with $\partial\Sigma=\Gamma$.	
\end{rmk}




\section{Counterexamples}

In this section, we discuss potential generalizations, and construct some important examples in these directions. First, we discuss the embeddedness of the solution to the Plateau problem for weak-extreme curves.

\subsection{Embeddedness of least area disks bounding Weak-Extreme Curves} \label{Weakcounterexamplesec} \

\vspace{.2cm}

In our main result, we proved that a weak-extreme curve bounds an embedded \textit{minimal disk}. However, this does not imply the embeddedness of \textit{the least area disk} bounding the same curve. So, the following question is natural:

\begin{question} Let $\Gamma$ be a weak-extreme curve in $\BR^3$. Is the solution to the Plateau problem for $\Gamma$ is embedded?
	
\end{question}
	
Unfortunately, the answer to this question is "No". The following counterexample explains the main issue.

\begin{figure}[b]
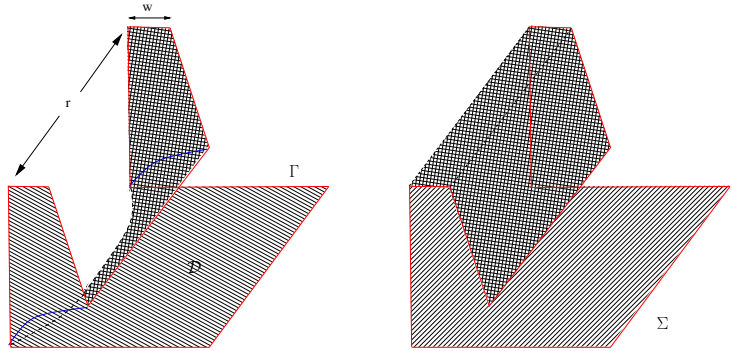

	\begin{center}
		
		$\begin{array}{c@{\hspace{.4in}}c}
		
		\resizebox{1.7in}{!}{\input{counterexample-1d.pstricks}} &
		
		\resizebox{1.7in}{!}{\input{counterexample-2c.pstricks}} \\
		
	\end{array}$
	
\end{center}

\caption{ \label{weakextremeCEfig} \footnotesize In the left, $\D$ is the least area disk in $\BR^3$ bounding the weak-extreme curve $\Gamma$ (red curve). In the right, $\Sigma$ is the stable minimal disk constructed in Theorem \ref{mainthm} bounding the same curve $\Gamma$.}
\end{figure}

In our construction, we decompose our curve $\Gamma$ into extreme pieces, and get embedded least area disk for each piece. By combining them, we got a piecewise smooth minimal disk $\wh{\Sigma}$. Then, we construct a special mean convex neighborhood near this piecewise smooth disk $\wh{\Sigma}$ by using area minimizing surfaces in both sides. Therefore, our smooth embedded minimal disk $\Sigma$ is basically "smoothening" of the piecewise smooth disk $\wh{\Sigma}$. 

However, the least area disk $\D$ in $\BR^3$ with $\partial \D=\Gamma$  might be completely different from our stable minimal disk $\Sigma$. In other words, while $\Sigma$ is the least area disk in a very special domain $\Z$, it may be far away from $\D$, the least area disk in whole $\BR^3$ (See Figure \ref{weakextremeCEfig}). As having smallest area is a somewhat weak condition to control the disk by the boundary curve, it is possible to construct a weak-extreme curve $\Gamma$ where the least area disk $\D$ in $\BR^3$ is very far from the stable embedded minimal disk $\Sigma$ constructed in Theorem \ref{mainthm}.

\vspace{.2cm}

\noindent {\bf Counterexample for Weak-Extreme Curves:} 

\vspace{.2cm}

In Figure \ref{weakextremeCEfig}, we have a weak-extreme curve $\Gamma$ (red curve). In the left, $\D$ is the least area disk in $\BR^3$ with $\partial \D=\Gamma$. In the right, $\Sigma$ is the stable minimal disk constructed in our main result (Theorem \ref{mainthm}). Notice that when $r$ is large, and $w$ is small, the least area disk takes the shape in the left.

\begin{figure}[h]
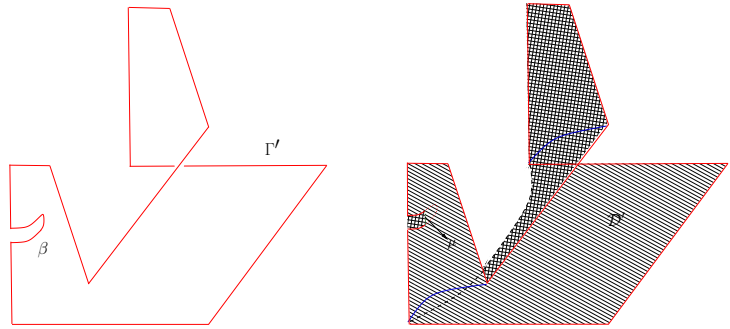

	\begin{center}
		$\begin{array}{c@{\hspace{.4in}}c}

		\resizebox{1.7in}{!}{\input{counterexample-3b.pstricks}}
		
		&
		
		\resizebox{1.7in}{!}{\input{counterexample-4b.pstricks}}
		
		\\
		
	\end{array}$
	
\end{center}

\caption{ \label{weakextremeCEfig2} \footnotesize In the left, we add a tiny hook to $\Gamma$ above, and obtain $\Gamma'$. While the stable minimal disk $\wh{\Sigma}'$ is still embedded for $\Gamma'$, the least area  disk $\D'$ is no longer embedded.}
\end{figure}

Now, by considering the difference between the stable minimal disk $\Sigma$ and the least area disk $\D$ above, we get the nonembedded least area disk for a weak-extreme curve as follows. Modify $\Gamma$ by adding a new hook as in the figure below. With the addition of this new hook $\beta$, we get a new weak-extreme curve $\Gamma'$ (See Figure \ref{weakextremeCEfig2}). $\Gamma'$ still bounds a stable minimal embedded disk which is a slight modification of $\Sigma$ above. However, the hook $\beta$ can be taken thin and long so that the least area disk $\D'$ is no longer embedded (In the figure right, $\D'$ has self-intersection along $\mu$).

In other words, the new weak-extreme curve still bounds a stable minimal disk $\Sigma'$ close to $\Sigma$ above by our main result. However, the least area disk $\D'$ in $\BR^3$ with $\partial \D'=\Gamma'$ is no longer embedded. This shows that being weak-extreme does not guarantee the embeddedness of the solution to the Plateau problem.

\subsection{Embeddedness of least area disks for Curves with Convex Hull Genus 0} \label{CHgenussec} \

\vspace{.2cm}

In this part, we  discuss a simple, but effective invariant of a Jordan curve to detect whether it bounds an embedded minimal disk: {\em Convex Hull Genus.}

\begin{defn} [Convex Hull Genus] For a given simple closed curve $\Gamma$ in $\BR^3$, let $\A_\Gamma$ be the space of embedded oriented surfaces in $CH(\Gamma)$ bounding $\Gamma$.  Let $g(\Gamma)=min\{genus(\Sigma)\mid \Sigma\in \A_\Gamma\}$. We call $g(\Gamma)$  {\em the convex hull genus of $\Gamma$}.		
\end{defn}

For example, any nontrivial knot in $\BR^3$, the convex hull genus is automatically positive. This is because it cannot bound an embedded disk as it is a nontrivial knot. Indeed, the convex hull genus is at least the genus of the knot (genus of the Seifert surface) in that case. In other words, geometric complexity of the curve (convex hull genus) is bounded below by its topological complexity (Seifert genus). There are also many examples of unknots in $\BR^3$ with positive convex hull genus. For further details on the convex hull genus, see \cite{Hu}.

Notice that for any minimal surface $S$ with $\partial S=\Gamma$, $S$ must belong to $CH(\Gamma)$ because of the maximum principle. Hence, if $g(\Gamma)>0$, then $\Gamma$ cannot bound any embedded minimal disk as there is no embedded disk in the convex hull. This means when convex hull genus is positive, it is a very useful invariant. By using this simple notion, Almgren and Thurston produced many examples of unknotted curves bounding high genus minimal surfaces \cite{AT}. 

We see that when convex hull genus is positive, the curve can not bound an embedded minimal disk. Hence, it is natural to ask how powerful this invariant is.

\begin{question} Let $\Gamma$ be a simple closed curve in $\BR^3$ with $g(\Gamma)=0$. Is the solution to the Plateau problem for $\Gamma$ embedded?	
\end{question}

Unfortunately, the answer to this question is also "No". The reason for that convex hull genus is indeed a very weak invariant to detect embeddedness. This is because one can enlarge the convex hull of a given curve by adding a very thin long tail. In particular, by adding such a very thin long tail to the curve, while one does not change the least area disk too much, the convex hull genus can easily be reduced to 0 for an unknot. The following example explains the situation.

\vspace{.2cm}

\noindent {\bf Counterexample for Convex Hull Genus $0$ Curves:}

\vspace{.2cm}

\begin{figure}[b]
	\begin{center}
		$\begin{array}{c@{\hspace{.3in}}c}

		\relabelbox  {\epsfxsize=2in \epsfbox{nonembedded.eps}}
		\relabel{1}{\small  $\Gamma$} 
		\relabel{2}{\small  $\Sigma$} 
		\relabel{3}{\small $\alpha$}
		
		\endrelabelbox &
		
		\relabelbox  {\epsfxsize=2in \epsfbox{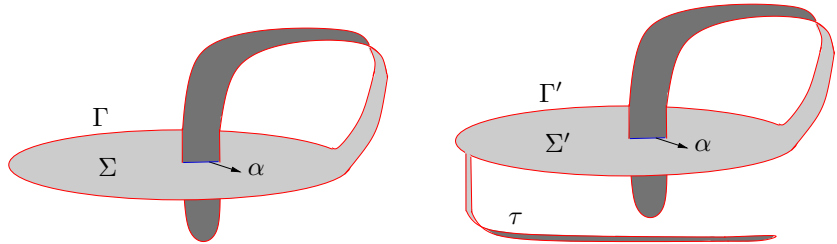}}
		\relabel{1}{\small  $\Gamma'$} 
		\relabel{2}{\small $\Sigma'$} 
		\relabel{3}{\small $\alpha$}
		\relabel{4}{\small $\tau$}
		\endrelabelbox \\		
		
	\end{array}$
	
\end{center}

\caption{ \label{weakextremeCEfig3} \footnotesize In the figure left, we decompose red curve $\Gamma$ into extreme pieces}
\end{figure}

Consider the nonembedded least area disk given in the figure below left. It is not hard to see that the convex hull genus of $\Gamma$ is positive. This is because lower tip of the hook is in the boundary of the convex hull, and any disk in the convex hull cannot go below that tip. 

On the other hand, convex hull genus can easily be manipulated by adding a thin tail to the curve like in the figure. By adding a very thin tail $\tau$ to $\Gamma$, we obtain a new Jordan curve $\Gamma'$  (See Figure \ref{weakextremeCEfig3}-right). As the extension $\tau$ is very thin, it slightly increases the area of the least area disk, i.e. $|\Sigma'|\sim |\Sigma|+ \epsilon$. However, the convex hull $CH(\Gamma')$ is much larger than the convex hull $CH(\Gamma)$. While the least area disks are basically same except the added tail, the convex hull genus of $\Gamma'$ is now $0$.

In particular, the solution to the Plateau problem  $\Sigma'$ for $\Gamma'$ is not embedded just like $\Sigma$. However, $g(\Gamma')=0$ as it is possible to find embedded disks bounding $\Gamma'$ in $CH(\Gamma')$, e.g. pull down the least area disk $\Sigma'$ in the figure along the self-intersection $\alpha$ until it is freed from the hook. Such an embedded disk cannot be in $CH(\Gamma)$ as the tip of the hook in the figure left belongs the the boundary of $CH(\Gamma)$. However, in the figure right, the tip of the hook is in the interior of the enlarged convex hull $CH(\Gamma')$, and hence it is possible to find embedded disks bounding $\Gamma'$ in $CH(\Gamma')$. In particular, $\Gamma'$ has convex hull genus $0$, but the solution to the Plateau problem is not embedded. Note also that $\Gamma'$ is not weak extreme either, as the hook stabs the core (Condition 3 violated).

\section{Concluding Remarks}

In this part, we  discuss further directions on the embeddedness problem.

\subsection{Difficulty of the Embeddedness Problem} \

\vspace{.2cm}

As mentioned above, so far there are only two strong result on the embeddedness of the solutions of the Plateau problem. The first one is Meeks-Yau's extreme curve result (Lemma \ref{MYlem}), and the second one is Jordan curves with total curvature less than $4\pi$ \cite{EWW}. Other than these results, there is not even a good conjecture on the problem.

Even though these two results are strong, they are far from being comprehensive. In particular, let $\D$ be an embedded least area disk in $\BR^3$. Let $\gamma$ be a simple closed curve in $\D$. Then, the solution to the Plateau problem for $\gamma$ would trivially be the subdisk $E\subset \D$ with $\partial E=\gamma$. Indeed, $E$ is the unique least area disk which $\gamma$ bounds by Meeks-Yau exchange roundoff trick. Since this is true for any such $\gamma\subset \D$, we can make $\gamma$ as irregular as we want while staying in $\D$. In particular, total curvature of $\gamma$ can be very large, and it may be far from being an extreme curve. Such examples show that the problem of embeddedness of the Plateau solution is still quite open.

So, the natural question is {\em "Why don't we even have a good conjecture on the embeddedness question?"}. The main difficulty basically lies in the Figure \ref{nonembeddedfig}. While a curve $\Gamma$ can bound an embedded least area disk $\Sigma$, if one adds a very thin small hook to $\Gamma$ which cuts through $\Sigma$, the new curve may no longer bounds an embedded minimal disk. Since the area does not give too much control on the curve, one can loose embeddedness with such simple operations. Hence, one needs to develop a notion where these two {\em seemingly close} curves are indeed distant. While adding a thin small hook  geometrically makes a major change in the curve, current geometric constraints (like $k(\Gamma)<4\pi$) are highly restrictive. 


\subsection{Generalization to $3$-Manifolds} \label{3mansec} \

\vspace{.2cm}

After getting the above result for $\BR^3$, it is natural to ask its generalization to general Riemannian $3$-manifolds. To use similar techniques in this setting, one needs the following two assumptions:

\vspace{.2cm}

\noindent {\em $\diamond$ $M$ is homogeneously regular:} This is a natural condition to make sure that the existence of least area disks \cite{HS}. $M$ is {\em homogeneously regular} means that $(M,g)$ has a upper bound for sectional curvature, and positive lower bound for the injectivity radius. Any compact manifold is automatically homogeneously regular. 
	
\vspace{.2cm}

\noindent {\em  $\diamond$ The convex hull of $\Gamma$ in $M$:} One can adapt the same definition of convex hull $CH(\Gamma)$ as the smallest convex set in $M$ containing $\Gamma$. Similarly, one can adapt extreme curve definition to this setting. However, the topology of $CH(\Gamma)$ may no longer be trivial in $M$ as in $\BR^3$. As we essentially use the triviality of the topology of $CH(\Gamma)$, to apply the techniques in Theorem \ref{mainthm}, one needs to make sure that $CH(\Gamma)$ has the required topological properties. For example, assuming the trivial second homology for $M$ ($H_2(M)=\{0\}$) could be a such good assumption.


\subsection{Further Questions} \

\vspace {.2cm}


\noindent {\em Generalization of Weak-Extreme Curves:} Considering the proof of the main theorem, a natural question would be \textit{"Which conditions for weak-extremeness are indeed necessary for such a result?} For example, it is clear that $\beta_i\cap\wh{\Sigma}=\emptyset$ is a crucial condition to implement the techniques in Theorem \ref{mainthm}. See also Figure \ref{nonembeddedfig}. In Remark \ref{rmk-Wecond}, we see that all the conditions are important in our construction. However, with different methods, some of those conditions might be removed.

\vspace {.15cm}

\noindent {\em New notions to detect non-embeddedness of the least area disks:} So far the only valuable notion to detect non-embeddedness of the solution to the Plateau problem is the convex hull genus (Section \ref{CHgenussec}). When $g(\Gamma)>0$, it automatically implies that any minimal disk bounding $\Gamma$ must be non-embedded. So, the following question is very natural: \textit{Are there other notions for Jordan curves in $\BR^3$ to detect non-embeddedness of the Plateau problem?}

\vspace {.2cm}

\noindent {\em Embedded Plateau Problem:} In \cite{Co2}, we discussed embeddedness of the solution to the Plateau problem from a different point of view. In particular, we showed that if $\Gamma$ is a simple closed curve bounding an embedded disk in a closed $3$-manifold $M$, then there exists a disk $\Sigma$ in $M$ with boundary $\Gamma$ such that $\Sigma$ minimizes the area among the embedded disks with boundary $\Gamma$. Moreover, $\Sigma$ is smooth, minimal and embedded everywhere except where the boundary $\Gamma$ meets the interior of $\Sigma$. The singularities occur in this construction are related to thin obstacle problem.

\end{document}